\documentclass{elsarticle}

\usepackage{amsmath}
\usepackage{amsfonts}
\usepackage{amssymb}

\newcommand{\abs}[1]{\ensuremath{\left| #1 \right| }}
\newcommand{\qPs}{$q$--Pochhammer symbol}
\newcommand{\bhs}{basic hypergeometric series}
\newcommand{\rhs}{right hand side}
\newcommand{\lhs}{left hand side}
\newcommand{\wps}{well--poised side}

\newcommand{\wrt}{with respect to}
\newcommand{\od}{one dimensional}

\newcommand{\BL}{Bailey Lemma}
\newcommand{\Bl}{Bailey lattice}
\newcommand{\Bc}{Bailey chain}
\newcommand{\Bp}{Bailey pair}
\newcommand{\Wt}{Watson transformation}

\newcommand{\RRis}{Rogers--Ramanujan identities}
\newcommand{\AGis}{Andrews--Gordon identities}
\newcommand{\Epnt}{Euler's Pentagonal Number Theorem}
\newcommand{\RSi}{Rogers--Selberg identity}

\newcommand{\Wdegf}{Weyl degree formula}

\newcommand{\hgr}{hyperoctahedral group}
\newcommand{\Js}{Jackson sum}
\newcommand{\ci}{cocycle identity}
\newcommand{\Jtpi}{Jacobi triple product identity} 
\newcommand{\tns}{\ensuremath{_{10}\varphi_9}}
\newcommand{\sfs}{\ensuremath{_6\varphi_5}}

\newcommand{\ml}{multilateralization lemma}
\newcommand{\Jc}{Jackson coefficients}
\newcommand{\eJc}{elliptic Jackson coefficients}

\newtheorem{thm}{Theorem}
\newtheorem{lem}[thm]{Lemma}
\newdefinition{defn}{Definition}
\newdefinition{rmk}{Remark}
\newproof{pf}{Proof}

\begin{document}

\begin{frontmatter}

\journal{}
\date{December 21, 2009}

\title{A Multilateral Bailey Lemma and Multiple \\ Andrews--Gordon Identities}  

\author{Hasan Coskun}
\ead{hasan\_coskun@tamu-commerce.edu}
\ead[url]{http://faculty.tamu-commerce.edu/hcoskun}
\address{Department of Mathematics, Texas A\&M
  University--Commerce, Binnion Hall, Room 314, Commerce, TX 75429}  

\begin{abstract}
A multilateral \BL\ is proved, and multiple analogues of the Rogers--Ramanujan identities and Euler's Pentagonal Theorem are constructed as applications. The extreme cases of the \AGis\ are also generalized using the multilateral \BL\ where their final form is written in terms of determinants of theta functions. 
\end{abstract}

\end{frontmatter}

\section{Introduction}
\label{section1}
Let $B_{k,i}(n)$ denote the number of partitions of $n$ of the form $(b_1\ldots b_s)$ where $b_j-b_{j+k-1}\geq 2$, and at most $i-1$ of $b_j$ equal 1. Let $A_{k,i}(n)$ denote the number of partitions of $n$ into parts that are not equivalent to $0, \pm i\, (\mathrm{mod} \,\,2k+1)$. B. Gordon~\cite{Gordon1} proved that $A_{k,i}(n)=B_{k,i}(n)$ for all $n$, and thus gave a combinatorial generalization for the \RRis.

G. E. Andrews~\cite{Andrews2} gave an analytic counterpart of Gordon's result, the Andrews--Gordon identities, and extended \RRis\ to all odd moduli. These identities may be written in the form
\begin{equation}
\label{AGis}
\sum_{n_1 \geq \cdots \geq n_{k-1} \geq 0} \dfrac{q^{n_1^2 +\cdots+
    n_{k-1}^2 + n_{i}+\cdots +n_{k-1} } } {(q)_{n_{1} -
    n_{2}} (q)_{n_{2} - n_{3}} \ldots (q)_{n_{k-1} } } 
= \prod_{\substack{n=1 \\ n\not\equiv 0, \, \pm i \, (\mathrm{mod} \,2k+1)}
    }^\infty \!\!\!\! \dfrac{1}{(1-q^n)}  
\end{equation}
where $i\in [k]$, $k\geq 2$ and $|q|<1$. Here, the \qPs\ $(a;q)_\alpha$ is defined formally by 
\begin{equation}
\label{qPochSymbol} 
(a;q)_\alpha :=\dfrac{(a;q)_\infty}{(aq^\alpha;q)_\infty}
\end{equation} 
where the parameters $a, q, \alpha\in\mathbb{C}$, and $(a;q)_\infty$ denotes the infinite product $(a;q)_\infty:=\prod_{i=0}^{\infty} (1-aq^i)$. 

The so--called extreme cases ($i=1$ and $i=k$) of the \AGis\
correspond to the $k$--th iteration of the one--parameter \BL. 
The full \AGis\ for any $i\in [k]$ may be proved by using the two--parameter \BL~\cite{AgarwalA1}. 

This paper is devoted to a multiple analogue of \AGis\ 
associated to the root systems of rank~$n$. We first construct a multilateral version of the $BC_n$ \BL~\cite{Coskun}, and then write a multiple analogue of \AGis. A brief review of the $BC_n$ \BL\ is presented in the next section. First, however, we recall the main properties of the well--poised Macdonald functions $W_{\lambda/\mu}(z; q,p,t; a,b)$ and Jackson coefficients $\omega_{\lambda/\mu}(z; r,q,p,t; a,b)$ used in the construction of the multiple \BL~\cite{Coskun1}. These functions are first introduced in the author's thesis~\cite{Coskun0} supervised by R. A. Gustafson.  

\section{Background}
\label{section5}

Let $V$ denote the 
space of infinite lower--triangular matrices whose entries are rational
functions 
over the field
$\mathbb{F}=\mathbb{C}(q, p, t,r,a,b)$ as in~\cite{Coskun1}. 
The condition that a matrix $u\in
V$ is lower triangular \wrt\ the partial inclusion
ordering $\subseteq$ defined by
\begin{equation}
\label{partialordering}
\mu \subseteq \lambda \;\Leftrightarrow \;\mu_i \leq \lambda_i, \quad
\forall i\geq 1.
\end{equation}
can be stated in the form
\begin{equation}
  u_{\lambda\mu} = 0 ,\, \quad \mathrm{when}\; \mu \not \subseteq
  \lambda.
\end{equation}
The multiplication operation in $V$ is defined by the relation
\begin{equation}
\label{multiplication}
  (uv)_{\lambda\mu} := \sum_{\mu\subseteq\nu\subseteq\lambda}
  u_{\lambda\nu} v_{\nu\mu}
\end{equation}
for $u,v\in V$. The matrices used in the definition of $BC_n$ \BL\ involve the elliptic
well--poised Macdonald functions $W_{\lambda/\mu}$ and elliptic well--poised \Jc\ $\omega_{\lambda/\mu}$
on $BC_n$. 

\subsection{Well--poised Macdonald functions}

Recall that an elliptic
analogue of the basic factorial is defined in~\cite{FrenkelT1} in the form 
\begin{equation}
  (a; q,p)_m := \prod_{k=0}^{m-1} \theta(aq^m)
\end{equation}
where $a\in \mathbb{C}$, $m$ is a positive integer and the normalized elliptic
function $\theta(x)$ is given by
\begin{equation}
  \theta(x) = \theta(x;p) := (x; p)_\infty (p/x; p)_\infty
\end{equation}
for $x, p\in \mathbb{C}$ with $\abs{p}<1$. The definition is extended to negative $m$ by setting $(a; q,p)_m = 1/ (aq^{m}; q, p)_{-m} $.  Note that when $p=0$, $(a; q,p)_m$ reduces to the standard (trigonometric) \qPs.

For any partition $\lambda = (\lambda_1, \ldots, \lambda_n)$ and
$t\in\mathbb{C}$, define~\cite{Warnaar2}
\begin{equation}
\label{ellipticQtPocSymbol}
  (a)_\lambda=(a; q, p, t)_\lambda := \prod_{k=1}^{n} (at^{1-i};
  q,p)_{\lambda_i} .
\end{equation}
Note that when $\lambda=(\lambda_1) = \lambda_1$ is a single part
partition, then $(a; q, p, t)_\lambda = (a; q, p)_{\lambda_1} =
(a)_{\lambda_1}$. The following notation will also be used.
\begin{equation}
  (a_1, \ldots, a_k)_\lambda = (a_1, \ldots, a_k; q, p, t)_\lambda :=
  (a_1)_\lambda \ldots (a_k)_\lambda .
\end{equation}

Now let $\lambda=(\lambda_1, \ldots, 
\lambda_n)$ and $\mu=(\mu_1, \ldots, \mu_n)$ be partitions of at most
$n$ parts for a positive integer $n$ such that the
skew partition $\lambda/\mu$ is a horizontal strip; i.e. $\lambda_1
\geq \mu_1 \geq\lambda_2 \geq \mu_2 \geq \ldots \lambda_n \geq
\mu_n \geq \lambda_{n+1} = \mu_{n+1} = 0$. Following~\cite{Coskun1}, define
\begin{multline}
\label{definitionHfactor}
H_{\lambda/\mu}(q,p,t,b) 
:= \prod_{1\leq i < j\leq
n}\left\{\dfrac{(q^{\mu_i-\mu_{j-1}}t^{j-i})_{\mu_{j-1}-\lambda_j}
(q^{\lambda_i+\lambda_j}t^{3-j-i}b)_{\mu_{j-1}-\lambda_j}}
{(q^{\mu_i-\mu_{j-1}+1}t^{j-i-1})_{\mu_{j-1}-\lambda_j}(q^{\lambda_i
    +\lambda_j+1}t^{2-j-i}b)_{\mu_{j-1}-\lambda_j}}\right.\\
\left.\cdot 
\dfrac{(q^{\lambda_i-\mu_{j-1}+1}t^{j-i-1})_{\mu_{j-1}-\lambda_j}}
{(q^{\lambda_i-\mu_{j-1}}t^{j-i})_{\mu_{j-1}-\lambda_j}}\right\}\cdot\prod_{1\leq
i <(j-1)\leq n}
\dfrac{(q^{\mu_i+\lambda_j+1}t^{1-j-i}b)_{\mu_{j-1}-\lambda_j}}
{(q^{\mu_i+\lambda_j}t^{2-j-i}b)_{\mu_{j-1}-\lambda_j}}
\end{multline}
and 
\begin{multline}
\label{definitionSkewW}
W_{\lambda/\mu}(x; q,p,t,a,b)
:= H_{\lambda/\mu}(q,p,t,b)\cdot\dfrac{(x^{-1}, ax)_\lambda
  (qbx/t, qb/(axt))_\mu}
{(x^{-1}, ax)_\mu (qbx, qb/(ax))_\lambda}\\
\cdot\prod_{i=1}^n\left\{\dfrac{\theta(bt^{1-2i}q^{2\mu_i})}{\theta(bt^{1-2i})}
  \dfrac{(bt^{1-2i})_{\mu_i+\lambda_{i+1}}}
{(bqt^{-2i})_{\mu_i+\lambda_{i+1}}}\cdot
t^{i(\mu_i-\lambda_{i+1})}\right\}.
\end{multline}
where $q,p,t,x,a,b\in\mathbb{C}$. The function
$W_{\lambda/\mu}(y, z_1, \ldots, 
z_\ell; q,p,t,a,b)$ is extended to $\ell+1$ variables $y, z_1, \ldots, z_\ell
\in\mathbb{C}$  
through the following recursion formula
\begin{multline}
\label{eqWrecurrence}
W_{\lambda/\mu}(y,z_1,z_2,\ldots,z_\ell;q, p, t, a, b) \\
= \sum_{\nu\prec \lambda} W_{\lambda/\nu}(yt^{-\ell};q, p, t, at^{2\ell},
bt^\ell) \, W_{\nu/\mu}(z_1,\ldots, z_\ell;q, p, t, a, b).
\end{multline}

We will also need the \eJc\ below. Let $\lambda$ and
$\mu$ be again 
partitions of at most $n$--parts such that $\lambda/\mu$ is a skew
partition. Then the \Jc\ $\omega_{\lambda/\mu}$ are defined by
\begin{multline}
\label{eq:omega{lambda,mu}}
\omega_{\lambda/\mu}(x; r, q,p,t; a,b)
:= \dfrac{(x^{-1}, ax)_{\lambda}} {(qbx, qb/ax)_{\lambda}}
    \dfrac{(qbr^{-1}x, qb/axr)_{\mu}}{(x^{-1}, ax)_{\mu}} \\
\cdot \dfrac{(r, br^{-1}t^{1-n})_{\mu}}{(qbr^{-2}, qt^{n-1})_{\mu}}
  \prod_{i=1}^{n}\left\{ \dfrac{\theta(br^{-1}t^{2-2i} q^{2\mu_i})}
    {\theta(br^{-1}t^{2-2i})}  \left(qt^{2i-2} \right)^{\mu_i} \right\} \\
\cdot \prod_{1\leq i< j \leq n} \hspace*{-5pt} \left\{ \dfrac{
    (qt^{j-i})_{\mu_i - \mu_j} } { (qt^{j-i-1})_{\mu_i - \mu_j} }
    \dfrac{ (br^{-1}t^{3-i-j})_{\mu_i + \mu_j} } {
    (br^{-1}t^{2-i-j})_{\mu_i + \mu_j} } \right\} \\
\cdot W_{\mu} (q^{\lambda}t^{\delta(n)}; q, p, t, bt^{2-2n}, br^{-1}t^{1-n})
\end{multline}
where $x,r,q,p,t, a,b\in\mathbb{C}$. From here and on, $\delta(n)$ denotes the $n$-tuple $\delta(n) := (n-1, \ldots, 1, 0)$ and $q^{\lambda}t^{\delta(n)}$ denotes $q^{\lambda}t^{\delta(n)} := (q^{\lambda_1}t^{n-1}, \,q^{\lambda_2}t^{n-2}, \ldots, q^{\lambda_n})$. 

Note that $W_{\lambda/\mu}(x; q,p, t, a,b)$ vanishes unless $\lambda/\mu$
is a horizontal strip. However, $\omega_{\lambda/\mu}(x; r;
a,b)=\omega_{\lambda/\mu}(x; r, q,p,t; a,b)$ may be non--zero when
$\lambda/\mu$ is not a horizontal strip. 

The operator characterization~\cite{Coskun1} of $\omega_{\lambda/\mu}$
yields a recursion formula for \Jc\ in the form 
\begin{equation}
\label{recurrence22}
  \omega_{\lambda/\tau}(y,z; r; a,b) := \sum_\mu
  \omega_{\lambda/\mu}(r^{-k}y; r; ar^{2k},
  b r^k ) \, \omega_{\mu/\tau}(z; r; a, b)
\end{equation}
where $y=(x_{1},\ldots, x_{n-k})\in\mathbb{C}^{n-k}$ and
$z=(x_{n-k+1},\ldots, x_n)\in\mathbb{C}^k$.

A key result used in the development of the $BC_n$ \BL, the \ci\ for
$\omega_{\lambda/\mu}$, is written in~\cite{Coskun1} in the form  
\begin{multline}
\label{cocycleIdentity}
  \omega_{\nu/\mu}((uv)^{-1};uv,q,p,t;a(uv)^2, buv) \\ = \sum_{\mu\subseteq
  \lambda \subseteq \nu} \omega_{\nu/\lambda}(v^{-1};v,q,p,t;a(vu)^2,bvu) \,
  \omega_{\lambda/\mu}(u^{-1};u,q,p,t;au^2,bu)   
\end{multline}
where the summation index $\lambda$ runs over partitions.

Using the recurrence relation~(\ref{recurrence22}) the
definition of $\omega_{\lambda/\mu}(x;r;a,b)$ can be extended 
from the single variable $x\in\mathbb{C}$ case to the multivariable
function $\omega_{\lambda/\mu}(z; r; a,b)$ with arbitrary number of variables
$z = (x_1,\ldots, x_n)\in\mathbb{C}^n$. That $\omega_{\lambda/\mu}(z;
r; a,b)$ is symmetric is also proved in~\cite{Coskun1} using a
remarkable elliptic $BC_n$ \tns\ transformation identity.

\subsection{Limiting Cases}

The limiting cases of the basic (the $p=0$ case of the elliptic) $W$ functions
$W_{\lambda/\mu}(x; q,t,a,b) = W_{\lambda/\mu}(x; q, 0,t,a,b)$ will be
used in computations in what follows. To simplify the exposition, some more notation would be helpful. Set 
\begin{equation}
  W^b_{\lambda/\mu}(x; q, t, a) := \lim_{b
  \rightarrow 0}\, W_{\lambda/\mu}(x; q, t, a, b)  
\end{equation}
and
\begin{equation}
  W^a_{\lambda/\mu}(x; q, t, b) := \lim_{a
  \rightarrow 0}\, a^{-|\lambda|+|\mu|}\, W_{\lambda/\mu}(x; q, t, a, b)  
\end{equation}
and, finally 
\begin{equation}
  W^{ab}_{\lambda/\mu}(x; q,t, s) :=
  \lim_{a\rightarrow 0} W_{\lambda/\mu}(x; q,t, a, as) 
\end{equation}
The existence of these limits can be seen from ($p=0$ case of) the
definition~(\ref{definitionSkewW}),  the recursion
formula~(\ref{eqWrecurrence}) and the limit rule
\begin{equation}
\label{LimitRule}
  \lim_{a\rightarrow 0}\, a^{|\mu|} (x/a)_{\mu} 
= (-1)^{|\mu|}\, x^{|\mu|} t^{-n(\mu)} q^{n(\mu')} 
\end{equation}
where $\abs{\mu}=\sum_{i=1}^n \mu_i$ and $n(\mu) =
\sum_{i=1}^n (i-1) \mu_i$, 
and $n(\mu') =\sum_{i=1}^n \binom{\mu_i}{2}$. These functions are closely related to the Macdonald polynomials~\cite{Macdonald3}, interpolation Macdonald polynomials~\cite{Okounkov1} and $BC_n$ abelian functions~\cite{Rains1}. 

We now make these definitions more precise. Let $H_{\lambda/\mu}(q,t,b)=H_{\lambda/\mu}(q,0,t,b)$, and for $x\in \mathbb{C}$ define
\begin{multline}
W^a_{\lambda/\mu}(x;q, t, b)
:=\lim_{a\rightarrow 0} a^{-|\lambda|+|\mu|} W_{\lambda/\mu}(x;q,
t, a, b) \\ 
= (-qb/x)^{-|\lambda|+|\mu|} q^{ -n(\lambda') + n(\mu') } 
H_{\lambda/\mu}(q,t,b) \hspace*{-6pt} \\
\cdot \dfrac{(x^{-1})_\lambda (qbx/t)_\mu}
{(x^{-1})_\mu (qbx )_\lambda} 
\prod_{i=1}^n\left\{\dfrac{(1-bt^{1-2i}q^{2\mu_i})}{(1-bt^{1-2i})}
  \dfrac{(bt^{1-2i})_{\mu_i+\lambda_{i+1}}}
{(bqt^{-2i})_{\mu_i+\lambda_{i+1}}} \right\}
\end{multline}
Using~(\ref{eqWrecurrence}) 
we get the following recurrence formula for $W^a_{\lambda/\mu}$ function
\begin{equation}
W^a_{\lambda/\mu}(y,z;q, t,b) \\
= \sum_{\nu\prec \lambda} t^{2\ell(|\lambda|-|\nu|)}
W^a_{\lambda/\nu}(yt^{-\ell};q, t, bt^{\ell}) \, 
W^a_{\nu/\mu}(z;q, t, b)
\end{equation}
for $y\in\mathbb{C}$ and $z\in\mathbb{C}^\ell$.  Similarly, let 
\begin{multline}\label{definitionHfactor}
H_{\lambda/\mu}(q,t) 
:= \lim_{b\rightarrow 0} H_{\lambda/\mu}(q,t,b) \\ 
= \prod_{1\leq i < j\leq
n}\left\{\dfrac{(q^{\mu_i-\mu_{j-1}}t^{j-i})_{\mu_{j-1}-\lambda_j} }
{(q^{\mu_i-\mu_{j-1}+1}t^{j-i-1})_{\mu_{j-1}-\lambda_j} } 
\dfrac{(q^{\lambda_i-\mu_{j-1}+1}t^{j-i-1})_{\mu_{j-1}-\lambda_j}}
{(q^{\lambda_i-\mu_{j-1}}t^{j-i})_{\mu_{j-1}-\lambda_j}}
\right\}
\end{multline}
Sending $b\rightarrow
0$ we define the $W^{b}_{\lambda/\mu}$ functions in the form
\begin{multline}\label{definitionSkewW}
W^{b}_{\lambda/\mu}(x; q, t, a)
:=\lim_{b\rightarrow 0} W_{\lambda/\mu}(x;q,t, a, b) \\ 
= t^{-n(\lambda)+|\mu|+n(\mu) } H_{\lambda/\mu}(q, t) \dfrac{(x^{-1}, ax)_\lambda } {(x^{-1}, ax)_\mu } \hspace*{92pt}
\end{multline}
The recurrence formula for $W^{b}_{\lambda/\mu}$ would then be
\begin{equation}
W^{b}_{\lambda/\mu}(y,z;q, t, a) \\
= \sum_{\nu\prec \lambda} 
W^{b}_{\lambda/\nu}(yt^{-\ell};q, t, at^{2\ell}) \, 
W^{b}_{\nu/\mu}(z;q, t, a)
\end{equation}
for $y\in\mathbb{C}$ and $z\in\mathbb{C}^\ell$. 

Finally, by setting $b=as$ in the definition of $W_{\lambda/\mu}$ function, and sending $a\rightarrow
0$ we define another family of symmetric rational functions $W^{ab}{\lambda/\mu}$ in the form
\begin{multline}
W^{ab}_{\lambda/\mu}(x; q, t, s)
:=\lim_{a\rightarrow 0} W_{\lambda/\mu}(x;q,t, a, as) \\ 
= t^{-n(\lambda)+|\mu|+n(\mu) } H_{\lambda/\mu}(q, t) \dfrac{(x^{-1})_\lambda (qs/(xt))_\mu}
{(x^{-1})_\mu (qs/x)_\lambda} \hspace*{59pt}
\end{multline}
Using~(\ref{eqWrecurrence}) again
we get the following recurrence formula for $W^{ab}_{\lambda/\mu}$ function
\begin{equation}
W^{ab}_{\lambda/\mu}(y,z;q, t, s) \\
= \sum_{\nu\prec \lambda} 
W^{ab}_{\lambda/\nu}(yt^{-\ell};q, t, st^{-\ell}) \, 
W^{ab}_{\nu/\mu}(z;q, t, s)
\end{equation}
where $y\in\mathbb{C}$ and $z\in\mathbb{C}^\ell$ as before. We will often use two more limiting cases $W^{s\uparrow}_{\lambda/\mu}$ and $W^{s\downarrow}_{\lambda/\mu}$ defined as follows. 
\begin{multline}
  W^{s\uparrow}_{\lambda/\mu}(x; q,t) :=
  \lim_{s\rightarrow \infty} s^{|\lambda|-|\mu|} W^{ab}_{\lambda/\mu}(x; q,t, s) \\
  =(-q/x)^{-|\lambda|+|\mu|} q^{ -n(\lambda') + n(\mu') }  H_{\lambda/\mu}(q, t) \dfrac{(x^{-1})_\lambda }{(x^{-1})_\mu } \hspace*{51pt}
\end{multline}
The recurrence formula for $W^{s\uparrow}_{\lambda/\mu}$ function turns out to be
\begin{equation}
\label{Wsuprec}
W^{s\uparrow}_{\lambda/\mu}(y,z;q, t) \\
= \sum_{\nu\prec \lambda} t^{\ell(|\lambda|-|\nu|)}
W^{s\uparrow}_{\lambda/\nu}(yt^{-\ell};q, t) \, 
W^{s\uparrow}_{\nu/\mu}(z;q, t)
\end{equation}
for $y\in\mathbb{C}$ and $z\in\mathbb{C}^\ell$. Similarly,
\begin{multline}
  W^{s\downarrow}_{\lambda/\mu}(x; q,t) :=
  \lim_{s\rightarrow 0} W^{ab}_{\lambda/\mu}(x; q,t, s) \\
= t^{-n(\lambda)+|\mu|+n(\mu) } H_{\lambda/\mu}(q, t) \dfrac{(x^{-1})_\lambda }
{(x^{-1})_\mu } \hspace*{107pt}
\end{multline}
The recurrence formula for $W^{s\downarrow}_{\lambda/\mu}$ function may be written as
\begin{equation}
W^{s\downarrow}_{\lambda/\mu}(y,z;q, t) \\
= \sum_{\nu\prec \lambda}
W^{s\downarrow}_{\lambda/\nu}(yt^{-\ell};q, t) \, 
W^{s\downarrow}_{\nu/\mu}(z;q, t)
\end{equation}
for $y\in\mathbb{C}$ and $z\in\mathbb{C}^\ell$. 

\subsection{\BL}
\BL\ was introduced first by W. N. Bailey~\cite{Bailey1} in 1944 as a
special case of a general series transformation argument, as he
successfully attempted to clarify the mechanism behind Rogers'
proof~\cite{Rogers1} of the famous Rogers--Ramanujan
identities. G. E. Andrews~\cite{Andrews2} gave a stronger version of
the Lemma emphasizing its iterative nature in 1984. P. Paule~\cite{Paule1}    
independently presented a bilateral version of the Lemma around the
same time.

A. K. Agarwal, G. E. Andrews and
D. M. Bressoud~\cite{AgarwalA1} introduced 
an abstract matrix formulation of the one--parameter \BL\ and extended
the notion of a \Bc\ to that of a \Bl. Consequently, this formulation 
gave rise to a two parameter \Bl\ in~\cite{Bressoud1}. 
The most recent \od\ extension of the \BL\ was
given by Andrews~\cite{Andrews2}, where he defined the concept of a
well--poised \Bc\ and a Bailey tree. 

Bressoud,
Ismail and Stanton~\cite{Bressoud1} established variants of \BL\ where they 
replaced the base $q$ to $q^k$, instead of changing the parameters during the iteration. 
These versions of the \BL\ proved to be very effective as
well. They also introduced a method of inserting linear factors into
the \Bc\ which eliminated the need for the \Bl\ in the
\BL\ proof of \AGis.  

The general matrix formulation of the \BL\ may be given as follows.
For $q,a\in\mathbb{C}$, the
infinite lower triangular matrices 
\[M(a)=\left(\dfrac{1}{(q)_{i-j} (aq)_{i+j}} \right)\]
and 
\[M^{-1}(a)=\left( (1-aq^i) \dfrac{ (-1)^{i-j}
  (aq)_{i+j-1} q^{\binom{i-j}{2}} }{(q)_{i-j}} \right) \]
are inverses of each other. This result is equivalent to a
terminating $_4\varphi_3$ summation theorem. Now if $\alpha$ and
$\beta$ form a Bailey pair, 
that is if $\beta=M(a)\,\alpha$, then the pair $\beta'$ and
$\alpha'$ defined by $\alpha'= S(a)\, \alpha$ and $\beta'=
N(a) \,\beta$ form a Bailey pair, where $S(a)$ is the
infinite diagonal matrix 
\[S(a) := \left( \dfrac{(\rho)_i (\sigma)_i
    (aq/\rho\sigma)^i}{(aq/\rho)_i (aq/\sigma)_i}\delta_{ij}
    \right)\] 
and
\[N(a) := \left( \dfrac{(aq/\rho\sigma)_{i-j} (\rho)_j (\sigma)_j
    (aq/\rho\sigma)^{j} } {(q)_{i-j} (aq/\rho)_i (aq/\sigma)_i}
    \right) \] 
where $\rho,\sigma\in\mathbb{C}$.

The limiting case as $\sigma, \rho \rightarrow \infty$ of this result, which we call terminating weak \BL, may be written~\cite{Coskun1} explicitely as 
\begin{equation} 
\label{eq:weakBL1}
  \sum_{m=0}^n \dfrac{b^m q^{m^2}}{(q)_{n-m} } \cdot \sum_{k=0}^m
  \dfrac{ 1}{(q)_{m-k} (qb)_{m+k}} \cdot \alpha_k
= \sum_{k=0}^n \dfrac{b^k q^{k^2}} {(q)_{n-k} (qb)_{n+k} } \cdot \alpha_k
\end{equation}
Sending, in addition, $n\rightarrow \infty$ gives~\cite{Andrews2} the non--terminating weak \BL\ as a series identity in the form  
\begin{equation}
\label{eq:weakBL2}
  \sum_{m=0}^\infty b^m q^{m^2} \cdot \sum_{k=0}^m
  \dfrac{ 1}{(q)_{m-k} (qb)_{m+k}} \cdot \alpha_k
= \dfrac{1}{(qb)_{\infty}} \sum_{k=0}^\infty b^k q^{k^2} \alpha_k
\end{equation}
for any Bailey pair $(\alpha, \beta)$. This is in one dimensional case the version given by Paule. Sending $n\rightarrow \infty$ gives bilateral non--terminating \BL\ as follows.
\begin{equation}
\label{eq:weakBL2}
  \sum_{m=0}^\infty b^m q^{m^2} \cdot \sum_{k=0}^m
  \dfrac{ 1}{(q)_{m-k} (qb)_{m+k}} \cdot \alpha_k
= \dfrac{1}{(qb)_{\infty}} \sum_{k=0}^\infty b^k q^{k^2} \alpha_k
\end{equation}

Using these limiting cases, Paule~\cite{Paule1} gave a bilateral version of the \BL, where he defines a \Bp\ $(\alpha, \beta)$ by  
\[\beta_n=\sum_{k=-\infty}^\infty \dfrac{\alpha_k } {(q)_{n-k}
    (q)_{n+k+\delta }} \]
for $\delta\in\{0,1 \}$, and shows that if
\begin{equation}
  \label{eq:bilateralBL}
  \beta'_n = \sum_{j=0}^\infty \dfrac{ q^{j^2 +\delta j} } {(q)_{n-j} }
     \beta_j 
\end{equation}
and $\alpha'_n = q^{k^2+\delta k} \alpha_n$, then $(\alpha', \beta')$ also form a \Bp. We give multiple analogues of these limiting versions of \BL\ in this paper. 

As pointed out above, the one--parameter \BL\ depends on a single parameter $b$
which remains unchanged throughout the chain. However, by replacing the
matrix $M(a)$ by the two--parameter matrix~\cite{Bressoud1} 
\[M(a,b)=\left(\dfrac{ (b/a)_{i-j} (b)_{i+j} (1-a q^{2j}) a^{i-j}}
  {(q)_{i-j} (aq)_{i+j} (1-a)} \right)\]
the \BL\ is extended to a two parameter identity. 
Iteration of this version, extending
the Bailey chain concept, generates the so--called Bailey lattice. The matrix $M(a,b)$ has interesting properties such as 
\[M(b,c)M(a,b)= M(a,c)\quad \mathrm{and}\quad M^{-1}(a,b)=M(b,a).\]
These properties follow as special cases from a powerful result
given by Bressoud~\cite{Bressoud1}. Namely, 
\[S^{-1}(c) \,M(c,d)\, S(c) = S^{-1}(b)\, M(b,d)\, S(b) \,M(c, b),\]
for $qbc = d \rho \sigma$. 
With this extension the two--parameter \BL\ states that, if the 
infinite sequences $\alpha$ and $\beta$ 
form a Bailey pair with respect to $M(a,b)$, that is if
$\beta=M(a,b)\,\alpha$, then the pair $\beta'$ and $\alpha'$ defined by
\[\beta'= S(c) S^{-1}(b) M(b,d) S(b)\, \beta \]
and
\[\alpha'= S(a) M(a,c)\,\alpha \]
also form a Bailey pair with respect to $M(c,d)$ provided that $qbc = d
\rho \sigma$.   

First multiple series generalizations of the Bailey Lemma were apparently given
by Milne and Lilly~\cite{LillyM3} for the root systems
$A_l$ and $C_l$ of rank $l$. They generalized matrix
formulation of the \BL\ in these cases 
and gave numerous applications of their generalizations. A
different $A_2$--type \BL\ was presented more recently by Andrews,
Schilling and Warnaar~\cite{AndrewsS1} where they used supernomial
coefficients to replace the \od\ $q$--binomial coefficient in
bilateral \BL.   

Gustafson also gave a multiple analogue of the \BL\ 
in terms of his rational Schur functions $S_\lambda$ in an unpublished manuscript.  

The \BL\ proof of the extreme cases of \AGis\ may be summarized as follows. Start with
the elementary non--trivial sequence $\beta_n=\delta_{0n}$, and
using the inverse matrix $M^{-1}(b)$ compute the 
corresponding $\alpha_n$. Iterate the \BL\ with this pair $k$ times to write the generalized \Wt. 
Send the parameters $n, \sigma_i, \rho_i \rightarrow
\infty$ in the \Wt\ to get the generalized \RSi. Specialize the parameter $b$ 
by setting $b=1$ and $b=q$, and use the \Jtpi\ 
\begin{equation}
\label{JtpiBackground}
 (q, q/z, z)_\infty := \sum_{m=-\infty}^\infty (-1)^m q^{\binom{m}{2}}
 z^m,  \quad \abs{q}<1
\end{equation}
to write the final form~(\ref{AGis}) of the identities.  For the full \AGis\ one also changes the parameter $b$ during the iteration.

\subsection{One parameter $BC_n$ \BL}
We now give a review of our multiple~\cite{Coskun} \BL\ starting with the definitions of the multiple analogues of $M(b)$ and $S(b)$ matrices.   
\begin{defn} 
Let $\lambda$ be a partition of at most $n$--parts and
$b\in\mathbb{C}$. Define 
\begin{multline} 
\label{oneVarBM}
M_{\nu\lambda}(b):= (-1)^{|\lambda|} q^{2|\lambda|+n(\lambda')} 
  t^{ n(\lambda)+2(1-n)|\lambda|}  \,b ^{|\lambda|} \dfrac{ (b  t^{1-n})_\lambda }{ (qt^{n-1})_\lambda} \\
\cdot \prod_{1\leq i<j \leq n} \left\{
 \dfrac{ (qt^{j-i})_{\lambda_i -\lambda_j} }{(qt^{j-i-1})_{\lambda_i-\lambda_j}   }
 \dfrac{ (b  t^{3-i-j})_{\lambda_i+\lambda_j}} { (b  t^{2-i-j})_{\lambda_i+\lambda_j} } 
 \right\}  W^a_\lambda(q^\nu t^{\delta(n)}; q, t, b  t^{1-n})
\end{multline}
and the infinite diagonal matrix $S(b)$ with diagonal entries 
\begin{equation}
\label{Smatrix}
 S_{\lambda}(b) :=  (qb /\sigma\rho)^{|\lambda| } \dfrac{(\sigma, \rho)_{\lambda}}
 { (qb/\sigma, qb/\rho)_{\lambda} }
\end{equation}
where $\abs{\lambda}=\sum_{i=1}^n \lambda_i$ and $n(\lambda) =
\sum_{i=1}^n (i-1) \lambda_i$, 
and $n(\lambda') =\sum_{i=1}^n \binom{\lambda_i}{2}$.
\end{defn}

The properties of these matrices are investigated in~\cite{Coskun1}. It is shown, for example, that $M(b)$ is lower triangular and is independent of different representations of $\lambda$ in any dimension $n$. It was also shown that $M(b)$ is invertible where the inverse $M^{-1}(b)$ is an infinite dimensional lower triangular matrix with entries 
\begin{multline}
M^{-1}_{\lambda\mu}(b)=  \dfrac{ q^{-|\lambda|+|\mu|} t^{2 n(\mu)} }{ (qb , qt^{n-1} )_\mu } 
    \prod_{i=1}^{n} \left\{ \dfrac{(1-b  t^{2-2i} q^{2\lambda_i})  }{ (1-b  t^{2-2i})  } \right\} 
\\ \cdot  \prod_{1\leq i<j \leq n} \left\{
 \dfrac{ (qt^{j-i})_{\mu_i -\mu_j} }{(qt^{j-i-1})_{\mu_i-\mu_j}   }
 \right\}   W^b_\mu(q^\lambda t^{\delta(n)}; q, t, b  t^{2-2n})  
\end{multline}
Next, we recall~\cite{Coskun} the notion of a \Bp.
\begin{defn}
The infinite sequences $\alpha$ and $\beta$ of rational functions
$\alpha_{\lambda}, \beta_\lambda$ over the field $\mathbb{C}(q, t,r,a,b)$  
form a Bailey pair relative to $b$ if they satisfy   
\begin{equation}
  \beta_{\lambda} = \sum_{\mu} M_{\lambda\mu}(b)\,
  \alpha_{\mu}  
\end{equation}
where the sum is over partitions.
\end{defn}

The one parameter $BC_n$ \BL\ given in~\cite{Coskun} may now be stated as follows. 
\begin{thm}
\label{OneParaBaileyLem}
Suppose that the infinite sequences $\alpha$ and $\beta$ form a Bailey
pair relative to $b$. Then $\alpha'$ and $\beta'$ also form
a Bailey pair relative to $b$ where  
\begin{equation}
  \alpha'_{\lambda} = S_{\lambda}(b)\, \alpha_{\lambda} 
\end{equation}
and
\begin{equation}
  \beta'_{\lambda} = \sum_{\mu} N_{\lambda\mu}(b)\, \beta_{\mu} 
\end{equation}
where the sum is over partitions, and the entries of the
matrix are given by 
\begin{multline} 
\label{N(b)matrix}
N_{\nu\mu}(b)= q^{|\mu|} t^{2 n(\mu)} \dfrac{(qb , qb / \sigma\rho)_\nu } { (qb /\sigma, qb /\rho)_\nu}  
  \dfrac{  (\sigma, \rho)_\mu  }{(qb , qt^{n-1})_\mu} \\
\cdot  \prod_{1\leq i<j \leq n} \left\{
 \dfrac{ (qt^{j-i})_{\mu_i -\mu_j} }{(qt^{j-i-1})_{\mu_i-\mu_j}   }
 \right\}  W^{ab}_\mu(q^\nu t^{\delta(n)}; q, t, \sigma\rho t^{n-1} /qb ) 
\end{multline}
\end{thm}

The power of \BL\ comes from its potential for indefinite
iteration. The lemma can be 
applied to a given \Bp\ $(\alpha, \beta)$ repeatedly producing an
infinite sequence of \Bp s
$(\alpha,\beta)\rightarrow(\alpha',\beta')\rightarrow  
(\alpha'',\beta'')\rightarrow\cdots$, what is called a \Bc. In fact,
a stronger result says that it is possible to walk along the \Bc\
in every direction as depicted in the following figure.  
\begin{figure}[h]
\label{BaileyChain}
\begin{displaymath}
\begin{array}{ccccccccccccc}
\cdots & \leftrightarrow & \alpha^{(-2)} & \leftrightarrow &
\alpha^{(-1)} & \leftrightarrow & \alpha & \leftrightarrow &
\alpha' & \leftrightarrow & \alpha'' & \leftrightarrow & \cdots \\
 & & \updownarrow & & \updownarrow & & \updownarrow & &
\updownarrow & & \updownarrow &  \\
\cdots & \leftrightarrow & \beta^{(-2)} & \leftrightarrow &
\beta^{(-1)} & \leftrightarrow & \beta & \leftrightarrow &
\beta' & \leftrightarrow & \beta'' & \leftrightarrow & \cdots
\end{array}
\end{displaymath}
\caption{\Bc}
\end{figure}

The lower triangular matrices $M(b)$, $N(b)$ and the diagonal
matrix $S(b)$, having no zero entries on their diagonal, are all
invertible. One move forward and backward in the first line in
Figure 1 by $S(b)$
and $S^{-1}(b)$, in the second line by $N(b)$ and $N^{-1}(b)$ and move
up and down by $M(b)$ and $M^{-1}(b)$. That $N(b)$ is
invertible follows immediately from the construction $N(b)=M(b)S(b)M^{-1}(b)$. 
Therefore, the 
entire \Bc\ is uniquely determined by a single node $\alpha^{(i)}$
or $\beta^{(i)}$ for any $i\in\mathbb{Z}$ in the chain.

This powerful iteration mechanism allows one to prove numerous
multiple \bhs\ and multiple $q$-series identities. 
We will write a multilateral version of the multiple \BL\ for the extreme cases of the \AGis.  

The \Bp\ $(\alpha, \beta)$ corresponding to the simplest non--trivial
sequence $\beta$ defined by $\beta_{\lambda} = \delta_{\lambda 0}$ 
is called the unit \Bp. 
The corresponding $\alpha$ sequence is easily be computed to be
\begin{equation}
\label{UBPalpha}
  \alpha_{\lambda} = \sum_{\mu} M^{-1}_{\lambda\mu}(b) \,\beta_{\mu} =
  q^{-|\lambda| } \prod_{i=1}^n
  \left\{\dfrac{(1-bt^{2-2i}q^{2\lambda_i})} {(1-b t^{2-2i})} \right\}
\end{equation}
using the the inverse matrix $M^{-1}(b)$.

Iterating the multiple \BL\
$N$ times, starting with the unit \Bp\ 
yields the generalized \Wt\ in the form
\begin{multline}
\label{generalWatson}
\dfrac{ (qb, qb/\rho_N\sigma_N)_{\mu^0} }
    {(qb/\sigma_N, qb/\rho_N)_{\mu^0} } \cdot \!\!\! \sum_{\mu^{N-1}
    \subseteq \cdots \subseteq  \mu^0} 
    \!\! \prod_{k=1}^{N-1} \left\{ q^{|\mu^k|} t^{2n(\mu^k)} \dfrac{
    (qb/\rho_{N-k} 
    \sigma_{N-k}, \sigma_{N-k+1}, \rho_{N-k+1})_{\mu^{k}} }
    {(qb/\sigma_{N-k}, qb/\rho_{N-k}, qt^{n-1})_{\mu^k} }  
 \right. \\ \left. \cdot \prod_{1\leq i<j \leq n} \left\{\dfrac{
    (qt^{j-i})_{\mu^k_i-\mu^k_j} } {(qt^{j-i-1})_{\mu^k_i-\mu^k_j} }
    \right\} \, W^{ab}_{\mu^k}(q^{\mu^{k-1}} t^{\delta(n)}; q,t, 
    \rho_{N-k+1}\sigma_{N-k+1}t^{n-1}/qb ) \right\} \\ =
    \sum_{\mu\subseteq \mu^0} 
    (-1)^{|\mu|} q^{|\mu|+n(\mu')} t^{n(\mu)} \dfrac{(bt^{1-n})_\mu }{
    (qt^{n-1})_\mu } \prod_{i=1}^n
    \left\{\dfrac{(1-bt^{2-2i}q^{2\mu_i})} {(1-b t^{2-2i})} \right\}
    \\ \cdot  \prod_{1\leq i<j \leq n} \left\{\dfrac{
    (qt^{j-i})_{\mu_i-\mu_j} (bt^{3-i-j})_{\mu_i+\mu_j}}
    {(qt^{j-i-1})_{\mu_i-\mu_j} (b t^{2-i-j})_{\mu_i+\mu_j} }
    \right\}\, W^a_\mu(q^{\mu^0} t^{\delta(n)}; q,t, bt^{1-n}) \\
    \cdot \prod_{k=1}^N \left\{ \dfrac{(\sigma_{N-k+1},
    \rho_{N-k+1})_{\mu} } {(qb/\sigma_{N-k+1}, qb/\rho_{N-k+1})_{\mu}
    } \left( \dfrac{qb}{\sigma_{N-k+1}\rho_{N-k+1}} \right)^{|\mu|}
    \right\} 
\end{multline}
where
\begin{equation}
  \sum_{{\nu_{n}}\subseteq{\nu_{n-1}}\subseteq \cdots \subseteq
  {\nu_1} \subseteq\nu_{0} } :=
\sum_{{\nu_1}\subseteq\nu_{0}} \sum_{{\nu_2}\subseteq{\nu_1}}  \cdots
  \sum_{{\nu_{n}}\subseteq{\nu_{n-1}}} 
\end{equation}
It was shown~\cite{Coskun} that the $N=1$ case of the generalized \Wt~(\ref{generalWatson})
reduces to a multiple analogue of the terminating \sfs\ summation formula, and the
$N=2$ case reduces to that of the \Wt. 

Sending $\sigma_i,\rho_i\rightarrow \infty$ and $\mu^0\rightarrow
\infty$ in the
identity~(\ref{generalWatson}) 
gives the generalized \RSi\ 
in the next result.  
\begin{lem}
With notation as above, we have
\begin{multline}
\label{generalRSi}
(qb)_{\infty^n} \sum_{ \mu^{N-1}
  \subseteq \cdots \subseteq \mu^1 , \, \ell(\mu^1)\leq n}
  \!\!\!\!\! q^{|\mu^{1}| + 2 n({\mu^{1}}') }
  t^{(1-n)|\mu^{1}| } b^{|\mu^{1}|}  \\ 
 \prod_{1\leq i<j \leq n} \left\{\dfrac{(t^{j-i+1})_{\mu^{1}_i 
-\mu^{1}_j} }{(t^{j-i})_{\mu^{1}_i -\mu^{1}_j} } \dfrac{
  (qt^{j-i})_{\mu^1_i-\mu^1_j} } {(qt^{j-i-1})_{\mu^1_i-\mu^1_j} }
  \right\} \dfrac{ 1 } { (qt^{n-1})_{\mu^1} }  \\ 
\prod_{k=2}^{N-1} \left\{ q^{|\mu^k|} t^{2n(\mu^k)}  
\prod_{1\leq i<j \leq n} \left\{\dfrac{ (qt^{j-i})_{\mu^k_i-\mu^k_j} }
{(qt^{j-i-1})_{\mu^k_i-\mu^k_j} } \right\} \dfrac{ 1 } {
  (qt^{n-1})_{\mu^k} } W^{s\uparrow}_{\mu^{k}}(q^{\mu^{k-1}} t^{\delta(n)}; q,t)
  \right\} \\  
= \sum_{\mu, \ell(\mu)\leq n} (-1)^{|\mu|} b^{N|\mu|}
  q^{N|\mu|+(2N+1) n(\mu')} t^{(-2N+1)n(\mu)+(1-n)|\mu|} \\
  \dfrac{(bt^{1-n})_\mu }{ (qt^{n-1})_\mu }  \prod_{i=1}^n
  \left\{\dfrac{(1-bt^{2-2i}q^{2\mu_i})} {(1-b t^{2-2i})} \right\} \\ 
\prod_{1\leq i<j \leq n} \left\{\dfrac{ (qt^{j-i})_{\mu_i-\mu_j}
    (bt^{3-i-j})_{\mu_i+\mu_j}} 
{(qt^{j-i-1})_{\mu_i-\mu_j} (b t^{2-i-j})_{\mu_i+\mu_j} }
  \dfrac{(t^{j-i+1})_{\mu_i -\mu_j} (qbt^{2-i-j})_{\mu_i+\mu_j}}
{(t^{j-i})_{\mu_i -\mu_j} (qbt^{1-i-j})_{\mu_i+\mu_j}} 
\right\} 
\end{multline}
\end{lem}

\begin{pf}
First we set $\mu^0 = M^n$ in the generalized
\Wt~(\ref{generalWatson}) and use the analogue of the
\Wdegf~(\ref{conj:degform})
and the limit rule~(\ref{LimitRule}), and send
$\sigma_i,\rho_i\rightarrow \infty$ to get 
\begin{multline}
(qb)_{M^n} \sum_{\mu^{N-1}
  \subseteq \cdots \subseteq \mu^1 \subseteq M^n} 
  \!\!\!\!\! (-1)^{|\mu^{1}|}
  q^{(M+1)|\mu^{1}| + n({\mu^{1}}') }
  t^{(1-n)|\mu^{1}|+ n(\mu^{1})} b^{|\mu^{1}|}   \\ 
 \prod_{1\leq i<j \leq n} \left\{\dfrac{(t^{j-i+1})_{\mu^{1}_i 
-\mu^{1}_j} }{(t^{j-i})_{\mu^{1}_i -\mu^{1}_j} } \dfrac{
  (qt^{j-i})_{\mu^1_i-\mu^1_j} } 
{(qt^{j-i-1})_{\mu^1_i-\mu^1_j} } \right\} \dfrac{ (q^{-M})_{\mu^1} } {
  (qt^{n-1})_{\mu^1} }  \\ 
\prod_{k=2}^{N-1} \left\{ q^{|\mu^k|} t^{2n(\mu^k)}  
\prod_{1\leq i<j \leq n} \left\{\dfrac{ (qt^{j-i})_{\mu^k_i-\mu^k_j} }
{(qt^{j-i-1})_{\mu^k_i-\mu^k_j} } \right\} \dfrac{ 1 } {
  (qt^{n-1})_{\mu^k} }   W_{\mu^{k}}(q^{\mu^{k-1}} t^{\delta(n)}; q,t)
  \right\} \\  
= \sum_{\mu\subseteq M^n}  (-1)^{(2N+2)|\mu|} b^{N|\mu|}
  q^{(N+M)|\mu|+2N n(\mu')} t^{(-2N+2)n(\mu)+(1-n)|\mu|}  \\
  \dfrac{(q^{-M})_\mu } {(q^{1+M}b)_\mu }  \dfrac{(bt^{1-n})_\mu }{
  (qt^{n-1})_\mu }  \prod_{i=1}^n
  \left\{\dfrac{(1-bt^{2-2i}q^{2\mu_i})} {(1-b t^{2-2i})} \right\} \\ 
\prod_{1\leq i<j \leq n} \left\{\dfrac{ (qt^{j-i})_{\mu_i-\mu_j}
    (bt^{3-i-j})_{\mu_i+\mu_j}} 
{(qt^{j-i-1})_{\mu_i-\mu_j} (b t^{2-i-j})_{\mu_i+\mu_j} }
  \dfrac{(t^{j-i+1})_{\mu_i -\mu_j} (qbt^{2-i-j})_{\mu_i+\mu_j}}
{(t^{j-i})_{\mu_i -\mu_j} (qbt^{1-i-j})_{\mu_i+\mu_j}} 
\right\} 
\end{multline}
We next send $M\rightarrow \infty$ as we apply the dominated
convergence theorem as stated in ~\cite{Coskun} to get 
the generalization of the $BC_n$ \RSi~(\ref{generalRSi}) to be
proved. The convergence theorem applies, because for $a,b\in\mathbb{C}$ such that the denominator never vanishes there exists a constant $C$ independent of $\alpha$ 
such that 
\begin{equation}
\abs{ \dfrac{(a q^\alpha)_\infty }{(b q^\alpha)_\infty} } < C    
\end{equation}
when $\abs{q}<1$. This observation also implies, in the view of the recurrence relation~(\ref{eqWrecurrence}), that 
for any partitions $\lambda, \mu$ of at most $n$--parts and $q, t\neq
0, a\neq 0, b \in\mathbb{C}$ such that $\abs{q}<1$, we have 
\begin{equation}
\abs{W_\mu (q^\lambda t^{\delta(n)}; q, t; a, b) } < C_w    
\end{equation}
where $C_w$ is a constant independent of $\lambda$ and $\mu$.
\end{pf}

Before writing the extreme cases of our multiple \AGis, we multilateralize the summand in the well--poised \rhs\
of~(\ref{generalRSi}).  
\begin{lem}
The summand in the well--poised \rhs\ of the generalized \RSi~(\ref{generalRSi}) 
is symmetric under the standard hyperoctahedral group action of permuattions and sign changes. 
\end{lem}

\begin{pf}
Set $q^{z_i} = b^{1/2} t^{i-1}$ in~(\ref{generalRSi}) and use the definition~(\ref{qPochSymbol}) 
and standard properties of the \qPs\ to write the summand in the form 
\begin{multline}
\label{genRSiHyperSym}
\prod_{i=1}^n q^{-(N-1)z_i^2 } \prod_{i=1}^{n}
  \dfrac{(q^{1+2z_i})_\infty (q^{1-2z_i})_\infty} 
  {(qt^{n-i})_\infty (q^{1-2z_i}t^{n-i})_\infty } \\  
  \prod_{1\leq i<j \leq n} \left\{ \dfrac{ (q^{1+z_i-z_j})_\infty
  (q^{1-z_i+z_j})_\infty } 
{(t^{-1}q^{1+z_i-z_j})_\infty (t^{-1}q^{1-z_i+z_j})_\infty } 
  \dfrac{ (q^{1-z_i-z_j})_\infty
  (q^{1+z_i+z_j})_\infty } 
{(t^{-1}q^{1-z_i-z_j})_\infty (t^{-1}q^{1+z_i+z_j})_\infty } \right\} \\
\sum_{\mu,\ell(\mu)\leq n} \prod_{i=1}^n
  q^{(N-1)(\mu_i+z_i)^2 }  \prod_{i=1}^{n}
  \dfrac{(q^{1-z_i+(\mu_i+z_i)}t^{n-i})_\infty
  (q^{1-z_i-(\mu_i+z_i)}t^{n-i})_\infty 
  } {(q^{1+2z_i+2\mu_i})_\infty (q^{1-2z_i-2\mu_i})_\infty}  \\  
  \prod_{1\leq i<j \leq n} \left\{\dfrac{
  (q^{1+z_i-z_j+\mu_i-\mu_j})_\infty 
  (q^{1-z_i+z_j-\mu_i +\mu_j})_\infty } 
{(t^{-1}q^{1+z_i-z_j+\mu_i-\mu_j})_\infty (t^{-1}q^{1-z_i+z_j-\mu_i
  +\mu_j})_\infty }  \right. \\ \left.
\dfrac{ (q^{1-z_i-z_j-\mu_i-\mu_j})_\infty 
  (q^{1+z_i+z_j+\mu_i +\mu_j})_\infty } 
{(t^{-1}q^{1-z_i-z_j-\mu_i-\mu_j})_\infty (t^{-1}q^{1+z_i+z_j+\mu_i
  +\mu_j})_\infty }  \right\}
\end{multline}
It is now obvious that the summand
has the \hgr\ symmetries $q^{z+\mu} \leftrightarrow q^{w(z+\mu)}$.  
\end{pf}

We now write a $D_n$ generalization of the \AGis~(\ref{AGis}) in
the two extreme cases ($i=1$ and $i=k$) in the next theorem.
\begin{thm}
\label{AGisExtremeCases}
Let $\abs{q}<1$ and $N$ be a positive integer. The \wps\ of the \AGis\
can be written in the form  
\begin{multline}
\label{AGis_z=n-i}
\prod_{i=1}^{n-1} \left\{\dfrac{1}{(1+q^{n-i})} \right\} 
\prod_{1\leq i<j \leq n} \left\{\dfrac{ 1} { (1-q^{j-i})^2
    (1-q^{2n-i-j})^2 } \right\} \dfrac{1}{2} \, (-1)^{\binom{n}{2}}
\prod_{i=1}^n q^{ (n-i)^2 } \\
 \det_{1\leq i,j\leq n} \bigg( q^{(j-1)(n-i)} (q^{2N+1}, q^{
     (2n-2i+1) N+(j-1)}, q^{2N+2 - (2n-2i+1) N-j};
   q^{2N+1})_\infty  \bigg. \\ \bigg. + q^{- (j - 1)(n-i ) }
(q^{2N+1}, q^{(2n-2i+1) N-(j-1)}, q^{2N - (2n-2i+1) N+j};
   q^{2N+1})_\infty  \bigg) 
\end{multline}
when $b=t^{2n-2}$.
For the specialization $b=qt^{2n-2}$, we get 
\begin{multline}
\label{AGis_z=n-i+1/2}
\prod_{i=1}^n \left\{\dfrac{q^{ (n-i)(n-i+1/2) }}
  {(1-q^{2n-2i+1} ) }  \right\} \prod_{1\leq i<j \leq n} \left\{
  \dfrac{1} {(1-q^{j-i} )^2 (1-q^{2n+1-i-j})^2 } 
\right\}  \dfrac{1}{2} \, (-1)^{\binom{n}{2}}  \\ 
  \det_{1\leq i,j\leq n} \bigg( q^{(j-1)(n-i+1/2)} (q^{2N+1},
  q^{(2n-2i+2) N + (j-1) }, q^{2N+2-(2n-2i+2) N - j } ;
  q^{2N+1} )_\infty \bigg. \\ \bigg. + q^{-(j - 1)(n-i+1/2)}
  (q^{2N+1}, q^{(2n-2i+2) N - (j-1) }, q^{2N- (2n-2i+2) N
  + j } ; q^{(2N+1)} )_\infty \bigg)
\end{multline}
\end{thm}

\begin{pf}
By a routine application of the \ml\ from~\cite{Coskun}, we can write the well--poised side of the~(\ref{generalRSi}) in the form
\begin{multline}
\label{AGis_z=n-iSpecialization}
\prod_{i=1}^{n-1} \left\{\dfrac{1}{(1+q^{n-i})} \right\} 
\prod_{1\leq i<j \leq n} \left\{\dfrac{ 1} { (1-q^{j-i})^2
    (1-q^{2n-i-j})^2 } \right\} \\
\sum_{\mu\in \mathbb{Z}^n} \prod_{i=1} (-1)^{\mu_i}
q^{(2Nn-N+1-n- (2N-1)(i-1)) \mu_i} q^{(2N+1) \binom{\mu_i}{2} }
\\ \prod_{1\leq i<j \leq n} \left\{ (1-q^{j-i+\mu_i-\mu_j})
    (1-q^{2n-i-j+\mu_i+\mu_j})
\right\} 
\end{multline}
when $b=t^{2n-2}$ corresponding to $z_i=n-i$ specialization.
Similarly, we write
\begin{multline}
\label{AGis_z=n-i+1/2Specialization}
\prod_{i=1}^n \left\{\dfrac{1} {(1-q^{2n-2i+1} ) }  \right\}
\prod_{1\leq i<j \leq n} \left\{ \dfrac{1}  
{(1-q^{j-i} )^2 (1-q^{2n+1-i-j})^2 } 
\right\} \\ \sum_{\mu, \ell(\mu)\leq n} \prod_{i=1}^n (-1)^{\mu_i}
    q^{\left( 2nN +(1-n)-(2N-1)(i-1) \right)\mu_i } q^{(2N+1)
    \binom{\mu_i}{2} } \\ 
\prod_{1\leq i<j \leq n} \left\{ (1-q^{j-i+\mu_i-\mu_j} )
    (1-q^{2n+1 -i-j+\mu_i+\mu_j}) \right\} 
\end{multline}
for $b=qt^{2n-2}$ corresponding to the $z_i=n-i+1/2$ specialization.
Next, we employ the determinant evaluations 
\begin{multline}
   \prod_{1\leq i<j \leq n} (1-x_i x_j^{-1}) (1-x_i x_j) \\ 
   = \dfrac{1}{2} \, (-1)^{\binom{n}{2}} \prod_{i=1}^n x_i^{n-i}  
   \det_{1\leq i,j\leq n} \left( x_i^{j-1} + x_i^{-(j-1)} \right) 
\end{multline}
for the root system $D_n$ of rank $n$, and apply the \Jtpi~(\ref{JtpiBackground}) to
write the specializations~(\ref{AGis_z=n-iSpecialization})
and~(\ref{AGis_z=n-i+1/2Specialization}) in the forms to be proved. Note that the balanced \lhs\ of the series can be put into a form where the sum runs over all $n$-tuples of non--negative integers. However, we will not pursue it here. 
\end{pf}

\section{Multilateral \BL}
I would like to present our multilateralization argument in one dimensional case to make the multiple analogue easier to read. The series formulation of the classical \BL\ may be written in explicit form as follows.
\begin{multline}
  \sum_{m=0}^n  \dfrac{(\sigma, \rho)_m (qb/\sigma\rho)_{n-m} (qb/\sigma\rho)^m} 
  { (q)_{n-m} (qb/\sigma, qb/\rho)_{n} }  
  \sum_{k=0}^m \dfrac{ 1}{(q)_{m-k} (qb)_{m+k}} \, \alpha_k  \\
= \sum_{k=0}^n \dfrac{1}{(q)_{n-k} (qb)_{n+k}} 
\dfrac{(\sigma)_k (\rho)_k (qb/\sigma\rho)^k } 
{(qb/\sigma)_{k} (qb/\rho)_{k} } \cdot \alpha_k
\end{multline}
Note that the identity may be writen as
\begin{multline}
 \sum_{m=0}^n  \dfrac{(\sigma, \rho)_m (q^{1+\delta}/\sigma\rho)_{n-m} (q^{1+\delta}/\sigma\rho)^m} 
  { (q)_{n-m} (q^{1+\delta}/\sigma, q^{1+\delta}/\rho)_{n} }  
 \sum_{k=-m-\delta}^m \dfrac{ 1}{(q)_{m-k} (q^{1+\delta})_{m+k}} \, \alpha_k  \\
=  \sum_{k=-n-\delta}^n \dfrac{1}{(q)_{n-k} (q^{1+\delta})_{n+k}} 
\dfrac{(\sigma)_k (\rho)_k (q^{1+\delta}/\sigma\rho)^{k} } 
{(q^{1+\delta}/\sigma)_{k} (q^{1+\delta}/\rho)_{k} } \cdot \alpha_k
\end{multline}
where $\delta\in\{0,1\}$. 
Since the maps $k\leftrightarrow -k-\delta$ generate the set of all integers, we write the sum over all integers. 

This is what we call `strong bilateral' \BL. If we send the parameters $\sigma,\rho\rightarrow \infty$, we get the weak version
\begin{multline}
 \sum_{m=0}^n  \dfrac{q^{m(m + \delta)  }} { (q)_{n-m}  }  
 \sum_{k=-m-\delta}^m \dfrac{ 1}{(q)_{m-k} (q^{1+\delta})_{m+k}} \, \alpha_k  
=  \sum_{k=-n-\delta}^n \dfrac{q^{k(k + \delta) }}{(q)_{n-k} (q^{1+\delta})_{n+k}} 
\,  \alpha_k
\end{multline}
If we also send $n\rightarrow \infty$, we get 
\begin{multline}
 \sum_{m=0}^\infty q^{m(m + \delta) }  
 \sum_{k=-m-\delta}^m \dfrac{ 1}{(q)_{m-k} (q^{1+\delta})_{m+k}} \, \alpha_k  
=  \dfrac{1 }{ (q^{1+\delta})_{\infty}}  \sum_{k=-\infty}^\infty q^{k(k + \delta)} 
 \, \alpha_k
\end{multline} 
which,under standard converging conditions, gives the non--terminating bilateral \BL\ listed above. It should be noted that this technique may be applied to write bilateral version for many well-poised hypergeometric series identities that satisfy the invariance property under the action of sign changes $z+\delta/2 \leftrightarrow w(z+\delta/2)$ where $w\in\mathbb{Z}_2$. We will illustrate this below for multiple analogues of the very--well poised \sfs\ and \Js\ identities. 

We now give a multilateral version of multiple \BL~\cite{Coskun}. 
It was already shown~\cite{Coskun} that the matrix entries $M_{\lambda\mu}(b)$ and $S_\lambda(b)$
are invariant under the \hgr\ action of permutations an sign changes when $\lambda$ is a rectangular partition $\lambda=k^n$. More precisely, it was shown that under the specialization $t=q^k$ and $b=q^{2z_i+2k(i-1)}$ where $z_i\in\mathbb{C}$ and $k\geq 0$ is a non--negative integer, the matrix entries are invariant under the action $(\mu_i+z_i)\leftrightarrow w(\mu_i+z_i)$ for all elements $w\in W$, the \hgr\ or rank $n$. It was further verified that this action generates the full weight lattice $\mathbb{Z}^n$ only if $z_i=m/2+k(n-i)$ for some non--negative integers $m, k\geq 0$. Here, we extend these results for an arbitarray partition $\lambda$. 
\begin{thm}
The specialized matrix entries $M_{\nu\lambda}(b)$ are invariant under the \hgr\ action of sign changes and permutations for the specializations $b=q^{m+2k(n-1)}$ and $t=q^k$ where $q\in\mathbb{C}$ and $m,k\geq 0$ when $\nu$ is an arbitrary partition. 
\end{thm}

\begin{pf}
It was shown~\cite{Coskun} that $W_\lambda$ functions are well--defined for any $\lambda\in\mathbb{Z}^n$ extending the original definition given for partitions. Therefore, we only need to verify the invariance for the $W_\lambda$ function that enters the definition of $M_{\nu\lambda}(b)$. 

The duality formula~\cite{Coskun} for $W_{\lambda}$ functions states that 
\begin{multline}
\label{eq:duality}
W_{\lambda}\left(k^{-1}q^\nu t^\delta;q,t,k^2a,kb\right)
\cdot\dfrac{(qbt^{n-1})_\lambda (qb/a)_\lambda}{(k)_\lambda
  (kat^{n-1})_\lambda} \\
\cdot \prod_{1\leq i < j\leq n} \left\{ \dfrac{(t^{j-i})_{\lambda_i
    -\lambda_j}(qa^{\prime}t^{2n-i-j-1})_{\lambda_i+\lambda_j}}
{(t^{j-i+1})_{\lambda_i
-\lambda_j}(qa^{\prime}t^{2n-i-j})_{\lambda_i+\lambda_j}} \right\}\\
=W_{\nu}\left(h^{-1}q^\lambda t^\delta;q, t,h^2a^{\prime},hb\right)
\cdot \dfrac{(qbt^{n-1})_\nu (qb/a^{\prime})_\nu}{(h)_\nu
  (ha^{\prime}t^{n-1})_\nu} \\
\cdot \prod_{1\leq i < j\leq n} \left\{ \dfrac{(t^{j-i})_{\nu_i
      -\nu_j}(qat^{2n-i-j-1})_{\nu_i+\nu_j}} {(t^{j-i+1})_{\nu_i
-\nu_j}(qat^{2n-i-j})_{\nu_i+\nu_j}} \right\}
\end{multline}
where $k=a^{\prime}t^{n-1}/b$ and $h = at^{n-1}/b$. Since $W_{\nu}$ on the right is a symmetric function, the \lhs\ is invariant under the permutations of $k^{-1}q^\nu_i t^{n-i}$, or that of $q^{m/2+k(n-i)+\nu_i}$ upon setting $k=-m/2$ and $t=q^k$. This is precisely what we need for the multilateralization of $BC_n$ \BL. Moreover, the identity~\cite{Coskun}
\begin{multline}
a^{|\lambda|} b^{-|\lambda|} q^{-|\lambda|}
t^{-n(\lambda)+(n-1)|\lambda|} W_{\lambda}(x_1^{-1}, \ldots, x_n^{-1}; q^{-1}, p, t^{-1}, a^{-1},
b^{-1})\\
= a^{-|\lambda|} b^{|\lambda|} q^{|\lambda|}
t^{n(\lambda)-(n-1)|\lambda|} W_{\lambda}(x_1, \ldots, x_n; q, p, t,
a, b)
\end{multline}
shows that the $W_{\lambda}$ is invariant under sign changes too. This is true, in particular, if we set $x_i=q^{\lambda_i} t^{n-i}$ as needed in \BL. The symmetries for the diagonal $S_\nu(b)$ is verified in~\cite{Coskun} for arbitrary partitions $\nu$.
\end{pf}

Now we give our multilateral definitions for the multiple $M(b)$ and $S(b)$ matrices.   
\begin{defn}
Let $\lambda$ be a partition of at most $n$--parts and
$b=q^{m+2k(n-1)}$ and $t=q^k$ for $q\in\mathbb{C}$ and $m,k\geq 0$. Define 
\begin{multline} 
\label{M(m,k)matrix}
M_{\nu\lambda}(m,k):=  q^{|\lambda|} 
  q^{ 2kn(\lambda)+ k(1-n)|\lambda| + m|\lambda|} \\
\cdot \!\!\! \prod_{1\leq i<j \leq n} \left\{
 \dfrac{ (q^{1+k(j-i)})_{\lambda_i -\lambda_j} }{(q^{1+k(j-i-1)})_{\lambda_i-\lambda_j}   }
 \dfrac{ (q^{m+k(2n+1-i-j)})_{\lambda_i+\lambda_j}} { (q^{m+k(2n-i-j)})_{\lambda_i+\lambda_j} } 
 \right\} \\ W^a_\lambda(q^\nu t^{\delta(n)}; q, q^k, q^{m+k(n-1)} )
\end{multline}
and
\begin{equation}
\label{S(m,k)matrix}
 S_{\lambda}(m,k) :=  \dfrac{ (q^{1+m+2k(n-1)} /\sigma\rho)^{|\lambda| } (\sigma, \rho)_{\lambda}} { (q^{1+m+2k(n-1)}/\sigma, q^{1+m+2k(n-1)}/\rho)_{\lambda} }
\end{equation}
where $\abs{\lambda}=\sum_{i=1}^n \lambda_i$ and $n(\lambda) =
\sum_{i=1}^n (i-1) \lambda_i$, 
and $n(\lambda') =\sum_{i=1}^n \binom{\lambda_i}{2}$ as before.
We also set
\begin{multline} 
\label{N(m,k)matrix}
N_{\nu\mu}(m,k) 
\\= \dfrac{(q^{1+m+2k(n-1)} , q^{1+m+2k(n-1)} / \sigma\rho)_\nu } 
{ (q^{1+m+2k(n-1)} /\sigma, q^{1+m+2k(n-1)} /\rho)_\nu}  
\dfrac{ q^{|\mu|+2k n(\mu)}  (\sigma, \rho)_\mu  }{(q^{1+m+2k(n-1)} , q^{1+k(n-1)})_\mu} \\
\cdot  \!\!\! \prod_{1\leq i<j \leq n} \left\{
 \dfrac{ (q^{1+k(j-i)})_{\mu_i -\mu_j} }{(q^{1+k(j-i-1)})_{\mu_i-\mu_j}   }
 \right\}  W^{ab}_\mu(q^\nu q^{k\,\delta(n)}; q, q^k, \sigma\rho q^{-1-m-k(n-1)} ) 
\end{multline}
\end{defn}
With these matrices, the strong multilateral \BL\ can be stated exactly as before. We will give the \BL\ in the special case when $m=\delta\in\{0,1\}$ as in the classical one dimensional case, and for $k=1$ or $t=q$. For clarity, we will drop $k$ from the notation when $k=1$ in the discussion below. $M_{\nu\lambda}(m)$, for example, denotes $M_{\nu\lambda}(m,1)$. In particular, for $k=1$ and $m=\delta$ we get
\begin{multline} 
M_{\nu\lambda}(\delta):= q^{ 2n(\lambda)+ (\delta+2-n)|\lambda| } \,
W^a_\lambda(q^\nu t^{\delta(n)}; q, q, q^{\delta+n-1} ) \\
\cdot \!\!\! \prod_{1\leq i<j \leq n} \left\{
 \dfrac{ (q^{1+j-i })_{\lambda_i -\lambda_j} }{(q^{j-i})_{\lambda_i-\lambda_j}   }
 \dfrac{ (q^{\delta+2n+1-i-j})_{\lambda_i+\lambda_j}} 
 { (q^{\delta+2n-i-j } )_{\lambda_i+\lambda_j} } 
 \right\}
\end{multline}
where $\delta\in\{0,1\}$. 

\begin{lem}
\label{multiOneParaBaileyLem}
Let the infinite sequences $\alpha=\{ \alpha_{\mu}: \mu\in\mathbb{Z}^n \}$ and $\beta=\{ \beta_\lambda:  \lambda \mathrm{\,\,is\, a\, partition\, of\, at\, most\,} n\mathrm{-part} \}$ of rational functions over the field $\mathbb{C}(q, t,r,a,b)$ 
form a Bailey pair relative to $b$. That is, they satisfy   
\begin{equation}
  \beta_{\lambda} = \sum_{\mu\in\mathbb{Z}^n} M_{\lambda\mu}(\delta)\,
  \alpha_{\mu}  
\end{equation}
where the sum terminates above at $\lambda$. Then $BC_n$ \BL\ implies that $\alpha'$ and $\beta'$ also form
a Bailey pair where  
\begin{equation}
\label{multiOneParaBaileyLem1}
  \alpha'_{\mu} = S_{\mu}(\delta)\, \alpha_{\mu} 
\end{equation}
for $\mu\in\mathbb{Z}^n$, and
\begin{equation}
\label{multiOneParaBaileyLem2}
  \beta'_{\nu} = \sum_{\lambda} N_{\nu\lambda}(\delta)\, \beta_{\lambda} 
\end{equation}
where the sum is over partitions. 
\end{lem}  
\begin{pf}
The proof is an immediate application of Lemma~\ref{OneParaBaileyLem}.
\end{pf}
This is our strong multilateral \BL. By sending $\sigma,\rho\rightarrow \infty$ and/or $\lambda\rightarrow \infty$ using the dominated convergence theorem both in $BC_n$ \BL\ of Theorem~\ref{OneParaBaileyLem} and in the multilateral $BC_n$ \BL\ of Lemma~\ref{multiOneParaBaileyLem}, we can write the terminating and non--terminating weak \BL s. We will only state multilateral terminating weak \BL\ here.

\begin{lem}
Let $M_{\nu\lambda}(\delta)$ 
matrix be defined as above in~(\ref{M(m,k)matrix}). Set
\begin{equation}
S_{\lambda}(\delta) = q^{(\delta+2(n-1))|\lambda|+n_2(\lambda)-2n(\lambda)} 
\end{equation}
where $n_2(\lambda) := |\lambda|+2n(\lambda') =
\sum_{i=1}^n\lambda_i^2$, and
\begin{multline} 
N_{\nu\mu}(\delta) = 
\dfrac{ q^{(\delta+n)|\mu| +n_2(\mu) } 
(q^{\delta+2n-1} )_\nu  } {(q^{\delta+2n-1},\, q^{n})_\mu } \\
\cdot \!\!\!
\prod_{1\leq i<j \leq n} \left\{
 \dfrac{ (1-q^{j-i+\mu_i -\mu_j} )}{(1-q^{j-i}) }
 \right\} W^{s\uparrow}_\mu(q^{\nu+\delta(n)}; q, q) 
\end{multline}
\end{lem}
The \BL\ of Lemma~\ref{multiOneParaBaileyLem} holds true with these definitions. 
\begin{pf}
The proof follows immediately from Theorem~\ref{multiOneParaBaileyLem} as $\sigma,\rho\rightarrow \infty$.
\end{pf}
Now, we write the inverse of $M(b)$ matrix.
\begin{lem}
For partitions $\lambda$ and $\mu$ of at most $n$-part, the inverse of $M(b)$ may  be written as 
\begin{multline}
M^{-1}_{\lambda\mu}(b):= 
(-1)^{|\lambda|}  t^{(n-1)|\lambda| -n(\lambda) } q^{n(\lambda')} 
 \dfrac{ (bt^{1-n})_\lambda }{
   (qt^{n-1})_\lambda}  
   \prod_{i=1}^n \left\{ \dfrac{ (1-b t^{2-2i} q^{2\lambda_i} )  }
   { (1-b t^{2-2i} )  } \right\}    \\ \cdot
  \dfrac{ q^{|\mu|} t^{2n(\mu)}  }{ (qt^{n-1})_\mu (qb)_\mu }    
 \prod_{1\leq i<j \leq n} \left\{
 \dfrac{ (qt^{j-i})_{\mu_i -\mu_j} }{(qt^{j-i-1})_{\mu_i-\mu_j}   }
  \right\}  W^b_\mu(q^\lambda t^{\delta(n)}; q, t, bt^{2-2n})  
\end{multline}
\end{lem}

\begin{pf}
This follows immediately form the \ci\ given in~\cite{Coskun}. 
\end{pf}

Note that under the specialization $b=q^{m+2k(n-1)}$ and $t=q^k$ for $m,k\geq 0$, we can write the inverse matrix as
\begin{multline}
M^{-1}_{\lambda\mu}(m,k)= 
(-1)^{|\lambda|}  q^{k(n-1)|\lambda| -kn(\lambda) + n(\lambda')} \\
\cdot \dfrac{ (q^{m+k(n-1)})_\lambda }{
   (q^{1+k(n-1)})_\lambda}  
   \prod_{i=1}^n \left\{ \dfrac{ (1- q^{m+2k(n-i)+2\lambda_i} )  }
   { (1-q^{m+2k(n-i)} )  } \right\}    
  \dfrac{ q^{|\mu|+2kn(\mu)}  }{ (q^{1+k(n-1)})_\mu (q^{1+m+2k(n-1)})_\mu }    \\
 \prod_{1\leq i<j \leq n} \left\{
 \dfrac{ (q^{1+k(j-i)})_{\mu_i -\mu_j} }{(q^{1+k(j-i-1)})_{\mu_i-\mu_j}   }
  \right\}  W^b_\mu(q^{\lambda+k\delta(n)}; q, q^k, q^m )  
\end{multline}
In particular, when $m=\delta$ and $k=1$, we get
\begin{multline}
M^{-1}_{\lambda\mu}(\delta)= 
(-1)^{|\lambda|}  q^{(n-1)|\lambda| -n(\lambda) + n(\lambda')} 
 \dfrac{ (q^{\delta+n-1})_\lambda }{
   (q^{n})_\lambda}  
   \prod_{i=1}^n \left\{ \dfrac{ (1- q^{\delta+2(n-i)+2\lambda_i} )  }
   { (1-q^{\delta+2(n-i)} )  } \right\}    \\ \cdot
  \dfrac{ q^{|\mu|+2n(\mu)}  }{ (q^{n})_\mu (q^{\delta+2n-1 })_\mu }    
 \prod_{1\leq i<j \leq n} \left\{
 \dfrac{ (q^{1+j-i})_{\mu_i -\mu_j} }{(q^{j-i})_{\mu_i-\mu_j}   }
  \right\}  W^b_\mu(q^{\lambda+ \delta(n)}; q, q, q^\delta )  
\end{multline}

With these definitions, we start iterating the multilateral \BL\ now. Note that if we iterate the strong multilateral \BL, we generate multiple multilateral versions of \bhs\ identities such as \sfs\ or \Js. However, we will present the weak versions here. 

The simplest non--trivial \Bp\ corresponds to $\beta_\lambda= \delta_{\lambda 0}$. In this case the corresponding $\alpha_\lambda$ sequence becomes
\begin{equation}
\alpha_\lambda=\sum_\mu M^{-1}_{\lambda\mu}(b) \beta_\mu  \\
  = (-1)^{|\lambda|}  t^{(n-1)|\lambda| -n(\lambda) } q^{n(\lambda')} \, f(\delta)
\end{equation}
where
\begin{equation}
f(\delta):= 
\dfrac{1}{n!} \prod_{i=1}^{n-1} \dfrac{ 1  }
   { (1+q^{n-i} )  }, \quad \mathrm{if\;} \delta=0 
\end{equation}
and
\begin{equation}
f(\delta):= 
\dfrac{1}{n!} \prod_{i=1}^{n} \dfrac{ 1  }
   { (1-q^{1+2n-2i} )  }, \quad \mathrm{if\;} \delta=1 
\end{equation}
Writing the first iteration 
explicitly gives 
\begin{multline} 
(q^{\delta+2n-1} )_\nu 
=  \sum_{\lambda\in\mathbb{Z}^n} 
W^a_\lambda(q^\nu t^{\delta(n)}; q, q, q^{\delta+n-1} ) \\
\cdot \!\!\! \prod_{1\leq i<j \leq n} \left\{
 \dfrac{ (q^{1+j-i })_{\lambda_i -\lambda_j} }{(q^{j-i})_{\lambda_i-\lambda_j}   }
 \dfrac{ (q^{\delta+2n+1-i-j})_{\lambda_i+\lambda_j}} 
 { (q^{\delta+2n-i-j } )_{\lambda_i+\lambda_j} } 
 \right\}  \\
q^{(2\delta+2n-1) |\lambda|} \,
 (-1)^{|\lambda|}  q^{n(\lambda')+n_2(\lambda)-n(\lambda)}  f(\delta)
\end{multline}
This identity yields the \Epnt\ in the limit as $\nu\rightarrow \infty$ when $\delta=0$. Using the identity~\cite{Coskun}
\begin{multline}
\lim_{k\rightarrow \infty} W^a_{\mu}(q^kt^{\delta(n)};q,q,q^{\delta+2n-2}q^{1-n}) \\
= (q^{\delta+n} )^{-|\mu|}   
\! \prod_{1\leq i < j\leq n} \dfrac{(q^{j-i+1})_{\mu_i -\mu_j}(q^{\delta+1+2n-i-j})_{\mu_i+\mu_j}}
{(q^{j-i})_{\mu_i -\mu_j} (q^{\delta+2n-i-j})_{\mu_i+\mu_j} } 
\end{multline}
and taking the limit, we get
\begin{multline} 
(q^{\delta+2n-1} )_{\infty^n}
=  \sum_{\lambda\in\mathbb{Z}^n} q^{(\delta+n-1) |\lambda|} \,
 (-1)^{|\lambda|}  q^{n(\lambda')+n_2(\lambda)-n(\lambda)}  \,n! f(\delta) \\
\cdot \!\!\! \prod_{1\leq i<j \leq n} \left\{
 \dfrac{ (1-q^{j-i+\lambda_i -\lambda_j} ) } {(1-q^{j-i})^2  }
 \dfrac{ (1-q^{\delta+2n-i-j+\lambda_i+\lambda_j}) } 
 { (1-q^{\delta+2n-i-j } )^2 } 
 \right\}
\end{multline}
which is the special $t=q$ case of the \Epnt\ given in~\cite{Coskun}. Iterate the \BL\ for a second time to get
\begin{multline} 
\sum_{\mu\subseteq \nu} \dfrac{ q^{(\delta+n)|\mu| +n_2(\mu) } 
(q^{\delta+2n-1} )_\nu  } {(q^{\delta+2n-1},\, q^{n})_\mu } \\
\cdot \!\!\!
\prod_{1\leq i<j \leq n} \left\{
 \dfrac{ (1-q^{j-i+\mu_i -\mu_j} )}{(1-q^{j-i}) }
 \right\} W^{s\uparrow}_\mu(q^{\nu+\delta(n)}; q, q) \cdot (q^{\delta+2n-1} )_\mu \\
=  \sum_{\lambda\in\mathbb{Z}^n} q^{ 2n(\lambda)+ (\delta+2-n)|\lambda| } \,
W^a_\lambda(q^{\nu+\delta(n)}; q, q, q^{\delta+n-1} ) \\
\cdot \!\!\! \prod_{1\leq i<j \leq n} \left\{
 \dfrac{ (q^{1+j-i })_{\lambda_i -\lambda_j} }{(q^{j-i})_{\lambda_i-\lambda_j}   }
 \dfrac{ (q^{\delta+2n+1-i-j})_{\lambda_i+\lambda_j}} 
 { (q^{\delta+2n-i-j } )_{\lambda_i+\lambda_j} } 
 \right\} \\ 
\cdot \left( q^{(\delta+2(n-1))|\lambda|+n_2(\lambda)-2n(\lambda)} \right)^2 
 (-1)^{|\lambda|}  q^{(n-1)|\lambda|-n(\lambda) } q^{n(\lambda')}  f(\delta)
\end{multline}
This is a multiple analogue of specialized \RSi~\cite{Coskun1}. Note also that although the series on the \rhs\ written over $\mathbb{Z}^n$, it actually terminates from above by $\nu$ and from below by 
$(-\nu_i - 2n- 2i+\delta) $.
In the limit $\nu\rightarrow \infty$, this identity gives the multiple analogues of  the celebrated first ($\delta=0$) and the second ($\delta=1$) \RRis. Recall the identity~\cite{Coskun} that
\begin{equation}
\lim_{k\rightarrow \infty} W^{s\uparrow}_{\mu}(q^kt^{\delta(n)};q,q) \\
= q^{-|\mu|} \! \prod_{1\leq i < j\leq n} \left\{ \dfrac{(q^{j-i+1})_{\mu_i
-\mu_j} } {(q^{j-i})_{\mu_i -\mu_j} } \right\}
\end{equation}
Therefore, in the limit we get
\begin{multline} 
\sum_{\mu\in\mathbb{Z}^n} \dfrac{ q^{(\delta+n-1)|\mu| +n_2(\mu) } 
 } {(q^{n})_\mu } \cdot \!\!\!
\prod_{1\leq i<j \leq n} \left\{
 \dfrac{ (1-q^{j-i+\mu_i -\mu_j} )}{(1-q^{j-i})^2 }
 \right\}  \\
= \dfrac{1}{(q^{\delta+2n-1} )_{\infty^n}} 
\sum_{\lambda\in\mathbb{Z}^n} \prod_{1\leq i<j \leq n} \left\{
 \dfrac{ (1-q^{j-i+\lambda_i -\lambda_j} )}{(1-q^{j-i})^2 } 
 \dfrac{ (1-q^{\delta+2n-i-j+\lambda_i+\lambda_j}) }
 { (1-q^{\delta+2n-i-j } )^2 } 
 \right\} \\ 
 \cdot \left( (-1)^{|\lambda|} q^{(2\delta+3(n-1))|\lambda|+2n_2(\lambda)-3n(\lambda)+n(\lambda')} \right) 
 n! f(\delta)
\end{multline}
This is precisely our multiple \RRis~\cite{Coskun1}. Repeating the iteration $N$ times in the same way generates the multiple \AGis\ given above in Theorem~\ref{AGisExtremeCases} for the extreme cases. 

\section{Conclusion}
The full version of the \AGis\ can be written in a similar way by using the mulatilateral version of the two-parameter \BL. We will write the full version of these identities in another publication.

\end{document}